\documentstyle[amscd,amssymb,verbatim,epsf]{amsart}
\begin{document}

\theoremstyle{plain}
\newtheorem{Thm}{Theorem}
\newtheorem{Cor}{Corollary}
\newtheorem{Con}{Conjecture}
\newtheorem{Main}{Main Theorem}

\newtheorem{Lem}{Lemma}
\newtheorem{Prop}{Proposition}

\theoremstyle{definition}
\newtheorem{Def}{Definition}
\newtheorem{Note}{Note}

\theoremstyle{remark}
\newtheorem{notation}{Notation}
\renewcommand{\thenotation}{}

\errorcontextlines=0
\numberwithin{equation}{section}
\renewcommand{\rm}{\normalshape}%

\title[On Variations of statistical ward continuity]
   {On Variations of statistical ward continuity}
\author{Huseyin Cakalli \\Maltepe University, Istanbul, Turkey}
\address{H\"Usey\.{I}n \c{C}akall\i\\
             Maltepe University, Department of Mathematics, Graduate School of Science and Engineering, Maltepe University, Marmara E\u{g}\.{I}t\.{I}m K\"oy\"u, TR 34857, Maltepe, \.{I}stanbul-Turkey \; \; \; \; \; Phone:(+90216)6261050 ext:2248, fax:(+90216)6261113} \email{hcakalli@@maltepe.edu.tr; hcakalli@@gmail.com}

\date{\today}

\begin{abstract}
In this paper, we introduce a concept of statistically $p$-quasi-Cauchyness of a real sequence in the sense that a sequence $(\alpha_{k})$ is statistically $p$-quasi-Cauchy if $\lim_{n\rightarrow\infty}\frac{1}{n}|\{k\leq n: |\alpha_{k+p}-\alpha_{k}|\geq{\varepsilon}\}|=0$ for each $\varepsilon>0$. A function $f$ is called statistically $p$-ward continuous on a subset $A$ of the set of real umbers $\mathbb{R}$ if it preserves statistically $p$-quasi-Cauchy sequences, i.e. the sequence $f(\textbf{x})=(f(\alpha_{n}))$ is statistically $p$-quasi-Cauchy whenever $\boldsymbol\alpha=(\alpha_{n})$ is a statistically $p$-quasi-Cauchy sequence of points in $A$. It turns out that a real valued function $f$ is uniformly continuous on a bounded subset $A$ of $\mathbb{R}$ if there exists a positive integer $p$ such that $f$ preserves statistically $p$-quasi-Cauchy sequences of points in $A$.
\end{abstract}
%%put here Key words

\maketitle

\section{Introduction}\label{sec1}
Throughout this paper, $\mathbb{N}$, and $\mathbb{R}$ will denote the set of positive integers, and the set of real numbers, respectively. $p$ will always be a fixed element of $\mathbb{N}$. The boldface letters such as $\boldsymbol\alpha$, $\boldsymbol\beta$, $\boldsymbol\zeta$ will be used for sequences $\boldsymbol\alpha=(\alpha_{n})$, $\boldsymbol\beta=(\beta_{n})$, $\boldsymbol\zeta=(\zeta_{n})$, ... of points in $\mathbb{R}$. A function $f:\mathbb{R} \longrightarrow \mathbb{R}$ is continuous if and only if it preserves convergent sequences. Using the idea of continuity of a real function in this manner, many kinds of continuities were introduced and investigated, not all but some of them we recall in the following: ward continuity (\cite{CakalliForwardcontinuity}, \cite{BurtonandColemanQuasiCauchysequences}), $p$-ward continuity (\cite{CakalliVariationsonquasiCauchysequences}),  $\delta$-ward continuity (\cite{CakalliDeltaquasiCauchysequences}), $\delta^{2}$-ward continuity (\cite{BrahaandCakalliAnewtypecontinuityforrealfunctions}), statistical ward continuity,  (\cite{CakalliStatisticalwardcontinuity}), $\lambda$-statistical ward continuity (\cite{CakalliandSonmezandAraslamdastatisticallywardcontinuity}), $\rho$-statistical ward continuity (\cite{CakalliAvariationonstatisticalwardcontinuity}), slowly oscillating continuity  (\cite{CakalliSlowlyoscillatingcontinuity,Vallin, CakalliandSonmezSlowlyoscillatingcontinuityinabstractmetricspaces}), quasi-slowly oscillating continuity (\cite{CanakandDik}), $\Delta$-quasi-slowly oscillating continuity (\cite{CakalliOnDeltaquasislowlyoscillatingsequences}), arithmetic continuity (\cite{YayingandHazarikaOnarithmeticcontinuity}, \cite{CakalliAvariationonarithmeticcontinuity}), upward and downward statistical continuities (\cite{CakalliUpwardanddownwardstatisticalcontinuities}), lacunary statistical ward continuity (\cite{CakalliandArasandSonmezLacunarystatisticalwardcontinuity}), $\delta$ lacunary statistical ward continuity (\cite{CakalliandKaplanAvariationonlacunarystatisticalquasiCauchysequencesCommunications}), $\delta^{2}$ lacunary statistical ward continuity (\cite{YildizIstatistikselboslukludelta2quasiCauchydizileri}), ideal ward continuity (\cite{CakalliandHazarikaIdealquasiCauchysequences, AkdumanandOzelandKilicmanSomeremarksoniwardcontinuity, CakalliAvariationonwardcontinuity}), ideal statistical ward continuity (\cite{SavasandCakalliIdealstatisticalquasiCauchysequences}), $N_{\theta}$-ward continuity (\cite{CakalliNthetawardcontinuity}, \cite{CakalliandKaplanAstudyonNthetaquasiCauchysequences}, \cite{CakalliandKaplanAvariationonstronglylacunarywardcontinuityJournalofMathematicalAnalysis}, \cite{KaplanandCakalliVariationsonstronglacunaryquasiCauchysequencesJNonlinearSciAppl}, \cite{KaplanandCakalliVariationsonstronglylacunaryquasiCauchysequencesAIP}), $N_{\theta}$-$\delta$-ward continuity (\cite{CakalliandKaplanAvariationonstronglylacunarywardcontinuityJournalofMathematicalAnalysis}), and Abel continuity (\cite{CakalliandAlbayrakNewtypecontinuitiesviaAbelconvergence}), which enabled some authors to obtain interesting results.

The purpose of this paper is to introduce statistically $p$-quasi-Cauchy sequences, and prove interesting theorems.

\section{Variations on statistical ward compactness}
The concept of a Cauchy sequence involves far more than that the distance between successive terms is tending to $0$ and statistically tending to zero, and more generally speaking, than that the distance between $p$-successive terms is statistically tending to zero, by $p$-successive terms we mean $\alpha_{k+p}$ and $\alpha_{k}$. Nevertheless, sequences which satisfy this weaker property are interesting in their own right.

Before giving our main definition we recall basic concepts. A sequence $(\alpha_{n})$ is called quasi Cauchy if $\lim_{n\rightarrow\infty}\Delta \alpha_{n}=0$, where $\Delta \alpha_{n}=\alpha_{n+1}-\alpha_{n}$ for each $n\in{\mathbb{N}}$ (\cite{BurtonandColemanQuasiCauchysequences}, \cite{CakalliForwardcontinuity}). The set of all bounded quasi-Cauchy sequences is a closed subspace of the space of all bounded sequences with respect to the norm defined for bounded sequences (\cite{NatarajanClassicalsummabilitytheory}). A sequence $(\alpha_{k})$ of points in $\mathbb{R}$ is slowly oscillating if $ \lim_{\lambda \rightarrow 1^{+}}\overline{\lim}_{n}\max _{n+1\leq k\leq [\lambda n]} |\alpha_{k}-\alpha_{n}| =0$, where $[\lambda n]$ denotes the integer part of $\lambda n$ (\cite{FDikMDikandCanak}). A sequence $(\alpha_{k})$ is quasi-slowly oscillating if $(\Delta \alpha_{k})$ is slowly oscillating. A sequence $(\alpha_{n})$ is called statistically convergent to a real number $L$ if $\lim_{n\rightarrow\infty}\frac{1}{n}|\{k\leq n: |\alpha_{k}-L|\geq{\varepsilon}\}|=0$ for each $\varepsilon>0$ (\cite{Fridy}, \cite{CakalliAstudyonstatisticalconvergence}, \cite{CakalliandKhanSummabilityintopologicalspaces}, and \cite{CasertaandKocinacOnstatisticalexhaustiveness}). Recently in \cite{CakalliVariationsonquasiCauchysequences} it was proved that a real valued function is uniformly continuous whenever it is $p$-ward continuous on a bounded subset of $\mathbb{R}$. Now we introduce the concept of a statistically $p$-quasi-Cauchy sequence.
\begin{Def} A sequence $(\alpha_{k})$ of points in $\mathbb{R}$ is called statistically $p$-quasi-Cauchy if $st-\lim_{k\rightarrow\infty} \Delta_{p} \alpha_{k}=0$, ie.  $\lim_{n\rightarrow\infty}\frac{1}{n}|\{k\leq n: |\Delta_{p} \alpha_{k}|\geq{\varepsilon}\}|=0$ for each $\varepsilon>0$, where $\Delta_{p} \alpha_{k}=\alpha_{k+p}-\alpha_{k}$ for every $k\in{\mathbb{N}}$.
\end{Def}
We will denote the set of all statistically $p$-quasi-Cauchy sequences by $\Delta^{s}_{p}$ for each $p\in{\mathbb{N}}$. The sum of two statistically $p$-quasi-Cauchy sequences is statistically $p$-quasi-Cauchy, the product of a statistically $p$-quasi-Cauchy sequence and a constant real number is statistically $p$-quasi-Cauchy, so that the set of all statistically $p$-quasi-Cauchy sequences  $\Delta^{s}_{p}$ is a vector space. We note that a  sequence is statistically quasi-Cauchy when $p=1$, i.e. statistically $1$-quasi-Cauchy sequences are statistical quasi-Cauchy sequences. It follows from the inclusion \\$\{k\leq n: |\alpha_{k+p}-\alpha_{k}| \geq\varepsilon\}\subseteq\\ \subseteq  \{k\leq n: |\alpha_{k+p}-\alpha_{k+p-1}| \geq \frac{\varepsilon}{p}\} \cup \{k\leq n: |\alpha_{k+p-1}-\alpha_{k+p-2}| \geq \frac{\varepsilon}{p}\}\cup ...\\ \cup \{k\leq n: |\alpha_{k+2}-\alpha_{k+1}| \geq \frac{\varepsilon}{p}\} \cup \{k\leq n: |\alpha_{k+1}-\alpha_{k}| \geq \frac{\varepsilon}{p}\}$\\
that any statistically quasi-Cauchy sequence is also statistically $p$-quasi-Cauchy, but the converse is not always true as it can be seen by considering the the sequence $(\alpha_{k})$ defined by $(\alpha_{k})=(0,1,0,1,...,0,1,...)$ is $2$-statistically quasi Cauchy which is not statistically quasi Cauchy. More examples can be seen  in \cite[Section 1.4]{NatarajanClassicalsummabilitytheory}. It is clear that any Cauchy sequence is in $\bigcap^{p=1}_{\infty}\Delta^{s}_{p}$, so that each $\Delta^{s}_{p}$ is a sequence space containing the space $\mathcal{C}$ of Cauchy sequences. It is also to be noted that $\mathcal{C}$ is a proper subset of $\Delta^{s}_{p}$ for each $p\in{\mathbb{N}}$.

\begin{Def} A subset $A$ of $\mathbb{R}$ is called statistically $p$-ward compact if any sequence of points in $A$ has a statistically $p$-quasi-Cauchy subsequence.
\end{Def}
We note that this definition of statistically $p$-ward compactness cannot be obtained by any summability matrix in the sense of \cite{CakalliSequentialdefinitionsofcompactness} (see also \cite{ConnorandGrosseSequentialdefinitionsofcontinuityforrealfunctions}, and \cite{CakalliOnGcontinuity}).

Since any statistically quasi-Cauchy sequence is statistically $p$-quasi-Cauchy we see that any ward compact subset of $\mathbb{R}$ is statistically $p$-ward compact for any $p\in{\mathbb{N}}$. A finite subset of $\mathbb{R}$ is statistically $p$-ward compact, the union of finite number of statistically $p$-ward compact subsets of $\mathbb{R}$ is statistically $p$-ward compact, and the intersection of any family of statistically $p$-ward compact subsets of $\mathbb{R}$ is statistically $p$-ward compact. Furthermore any subset of a statistically $p$-ward compact set of $\mathbb{R}$ is statistically $p$-ward compact and any bounded subset of $\mathbb{R}$ is statistically $p$-ward compact. These observations above suggest to us the following.
\begin{Thm} \label{ThmEisboundediffitispwardcompact}  A subset $A$ of $\mathbb{R}$ is bounded if and only if there exists a $p\in{\mathbb{N}}$ such that $A$ is statistically $p$-ward compact.
\end{Thm}
\begin{pf}
The bounded subsets of $\mathbb{R}$ are statistically $p$-ward compact, since any bounded sequence of points in a bounded subset of $\mathbb{R}$ is bounded and  any bounded sequence has a convergent subsequence which is statistically $p$-quasi-Cauchy for any $p\in{\mathbb{N}}$. To prove the converse, suppose that $A$ is not bounded. If it is unbounded above, pick an element $\alpha_{1}$ of $A$ greater than $p$. Then we can find an element $\alpha_{2}$ of $A$ such that $\alpha_{2}>2p+\alpha_{1}$. Similarly, choose an element $\alpha_{3}$ of $A$ such that $\alpha_{3}>3p+\alpha_{2}$. So we can construct a sequence $(\alpha_{j})$ of numbers in  $A$ such that  $\alpha_{j+1}>(j+1)p+\alpha _{j}$ for each $j\in{\mathbb{N}}$. Let $(\alpha _{j_{k}})$ be any subsequence of $(\alpha _{j})$. Since $\{k\leq n: |\alpha _{j_{k+1}}-\alpha _{j_{k}}|\geq\frac{p}{2}\}=\{1, 2, ..., n\}$, it follows that
$\lim_{n\rightarrow\infty}\frac{1}{n}|\{k\leq n: |\alpha _{j_{k+1}}-\alpha _{j_{k}}|\geq\frac{p}{2}\}|=1\neq{0}$. This means that $(\alpha _{j_{k}})$ is not statistically $p$-quasi-Cauchy. Then the sequence $(\alpha_{j})$ does not have any statistically $p$-quasi-Cauchy subsequence. If $A$ is bounded above and unbounded below, then pick an element $\beta_{1}$ of $A$ less than $-p$. Then we can find an element $\beta_{2}$ of $A$ such that $\beta_{2}<-2p+\beta_{1}$. Similarly, choose an element $\beta_{3}$ of $A$ such that $\beta_{3}<-3p+\beta_{2}$. Thus one can construct a sequence $(\beta_{i})$ of points in $A$ such that  $\beta_{i+1}<-(i+1)p+\beta _{i}$ for each $i\in{\mathbb{N}}$. Let $(\beta _{i_{k}})$ be any subsequence of $(\beta _{i})$. Since $\{k\leq n: |\beta _{i_{k+1}}-\beta _{i_{k}}|\geq p\}=\{1, 2, ..., n\}$, it follows that
$\lim_{n\rightarrow\infty}\frac{1}{n}|\{k\leq n: |\beta _{i_{k+1}}-\beta _{i_{k}}|\geq p\}|=1\neq{0}$. This implies that $(\beta _{i_{k}})$ is not statistically $p$-quasi-Cauchy. Then the sequence $(\alpha_{i})$ does not have any statistically $p$-quasi-Cauchy subsequence.
Thus this contradiction completes the proof of the theorem.
\end{pf}
It follows from Theorem \ref{ThmEisboundediffitispwardcompact} that statistically $p$-ward compactness of a subset of $A$ of $\mathbb{R}$ coincides with either of the following kinds of compactness: $p$-ward compactness (\cite[Theorem 2.3]{CakalliVariationsonquasiCauchysequences}), statistical ward compactness (\cite[Lemma 2]{CakalliStatisticalwardcontinuity}), $\lambda$-statistical ward compactness (\cite[Theorem 1]{CakalliandSonmezandAraslamdastatisticallywardcontinuity}), $\rho$-statistical ward compactness (\cite[Theorem 1]{CakalliAvariationonstatisticalwardcontinuity}), strongly lacunary ward compactness (\cite[Theorem 3.3]{CakalliNthetawardcontinuity}), slowly oscillating compactness (\cite[Theorem 3]{CakalliStatisticalquasiCauchysequences}), lacunary statistical ward compactness (see \cite{CakalliandArasandSonmezLacunarystatisticalwardcontinuity}), and \cite[Theorem 3]{CakalliStatisticalquasiCauchysequences}), ideal ward compactness (\cite[Theorem 8]{CakalliandHazarikaIdealquasiCauchysequences}), Abel ward compactness (\cite[Theorem 5]{CakalliandAlbayrakNewtypecontinuitiesviaAbelconvergence}).

If a closed subset of $\mathbb{R}$ is statistically $p$-ward compact for a positive integer $p$, then any sequence of points in $A$ has a $(P_n ,s)$-absolutely almost convergent subsequence (see \cite{CakalliGCanakPnsabsolutealmostconvergentsequences}, \cite{CakalliandTaylanOnabsolutelyalmostconvergence}, \cite{OzarslanandYildizAnewstudyontheabsolutesummabilityfactorsofFourierseries}, \cite{YildizAnewtheoremonlocalpropertiesoffactoredFourierseries}, \cite{BorOnGeneralizedAbsoluteCesaroSummability}, \cite{YildizOnapplicationofmatrixsummabilitytoFourierseries}, and \cite{YildizOnAbsoluteMatrixSummabilityFactorsofInfiniteSeriesandFourierSeries}).

\begin{Cor}
\normalfont
A subset of $\mathbb{R}$ is statistically $p$ ward compact if and only if it is statistically $q$ ward compact for any $p, q\in {\mathbb{N}}$.
\end{Cor}

\begin{Cor}
\normalfont
A subset of $\mathbb{R}$ is statistically $p$ ward compact if and only if it is both statistically upward half compact and statistically downward half compact.
\end{Cor}
\begin{pf}
The proof follows from \cite[Corollary 3.9]{CakalliUpwardanddownwardstatisticalcontinuities}, so is omitted.
\end{pf}

\begin{Cor}
\normalfont
A subset of $\mathbb{R}$ is  statistically  $p$ ward compact for a $p\in{\mathbb{N}}$ if and only if it is both lacunary statistically upward half compact and lacunary statistically downward half compact.
\end{Cor}
\begin{pf}
The proof follows from \cite[Theorem 1.3 and Theorem 1.9]{CakalliandMucukLacunarystatisticallyupwardanddownwardhalfquasiCauchysequencesJofAnalysis}, so is omitted.
\end{pf}

\section{Variations on statistical ward continuity}
 In this section, we investigate connections between uniformly continuous functions and statistically $p$-ward continuous functions. A function $f:\mathbb{R} \longrightarrow \mathbb{R}$ is continuous if and only if it preserves statistically convergent sequences. Using this idea, we introduce statistical  $p$-ward continuity.
\begin{Def}
A function $f$ is called statistically $p$-ward continuous on a subset $A$ of $\mathbb{R}$ if it preserves statistically $p$-quasi-Cauchy sequences, i.e. the sequence $f(\textbf{x})=(f(\alpha_{n}))$ is statistically $p$-quasi-Cauchy whenever $\boldsymbol\alpha=(\alpha_{n})$ is statistically $p$-quasi-Cauchy of points in $A$.
\end{Def}
We see that this definition of statistically $p$-ward continuity can not be obtained by any summability matrix $A$ (see \cite{ConnorandGrosseSequentialdefinitionsofcontinuityforrealfunctions}).

We note that the sum of two statistically $p$-ward continuous functions is statistically $p$-ward continuous, and for any constant $c\in \mathbb{R}$, $cf$ is statistically $p$-ward continuous whenever $f$ is a statistically $p$-ward continuous function, so that the set of all statistically $p$ ward continuous functions is a vector space. The composite of two statistically $p$-ward continuous functions is statistically $p$-ward continuous, but the product of two statistically $p$-ward continuous functions need not be statistically $p$-ward continuous as it can be seen by considering product of the statistically $p$-ward continuous function $f(x)=x$ with itself. If $f$ is a statistically $p$-ward continuous function, then $|f|$ is also statistically $p$-ward continuous since
$$|\{k\leq n:  |f(\alpha_{k+p})-f(\alpha_{k})|\geq{\varepsilon}\}|\subseteq |\{k\leq n: ||f(\alpha_{k+p})|-|f(\alpha_{k})||\geq{\varepsilon}\}|$$
which follows from the inequality $||f(\alpha_{k+p})|-|f(\alpha_{k})|| \leq  |f(\alpha_{k+p})-f(\alpha_{k})|$.
%for every $n\in {\mathbb{N}}$.
If $f$, and $g$ are statistically $p$-ward continuous, then $max\{f, g\}$ is also statistically $p$-ward continuous, which follows from the equality $max\{f, g\}=\frac{1}{2}\{|f-g|+|f+g|\}$.

In connection with statistically $p$-quasi-Cauchy sequences, slowly oscillating sequences, and convergent sequences the problem arises to investigate the following types of  "continuity" of a function on $\mathbb{R}$.
\begin{description}
\item[($\Delta^{s}_{p} $)] $(\alpha_{n}) \in {\Delta^{s}_{p}} \Rightarrow (f(\alpha_{n})) \in {\Delta^{s}_{p}}$
\item[($\Delta^{s}_{p} c$)] $(\alpha_{n}) \in {\Delta^{s}_{p}} \Rightarrow (f(\alpha_{n})) \in {c}$
\item[$(\Delta^{s})$] $(\alpha_{n}) \in {\Delta^{s}} \Rightarrow (f(\alpha_{n})) \in {\Delta^{s}}$
\item[$(c)$] $(\alpha_{n}) \in {c} \Rightarrow (f(\alpha_{n})) \in {c}$
\item[$(d)$] $(\alpha_{n}) \in {c} \Rightarrow (f(\alpha_{n})) \in {\Delta^{s}_{p}}$
\item[$(e)$] $(\alpha_{n}) \in {w} \Rightarrow (f(\alpha_{n})) \in {\Delta^{s}_{p}}$
\end{description}
where $w$ denotes the set of slowly oscillating sequences, and $\Delta^{s}=\Delta^{s}_{1}$. We see that $(\Delta^{s}_{p})$ is statistically $p$-ward continuity of $f$, and $(c)$ states the ordinary continuity of $f$. It is easy to see that $(\Delta^{s}_{p} c)$ implies $(\Delta^{s}_{p})$, and $(\Delta^{s}_{p} )$ does not imply $(\Delta^{s}_{p} c)$;  $\Delta^{s}_{p}$ implies $(d)$, and $(d)$ does not imply $(\Delta^{s}_{p})$;  and $(\Delta^{s}_{p})$ implies $(e)$, and $(e)$ does not imply $(\Delta^{s}_{p})$; $(\Delta^{s}_{p}c)$ implies $(c)$ and $(c)$ does not imply $(\Delta^{s}_{p} c)$; and $(c)$ implies $(d)$.
Now we see that $(\Delta^{s}_{p})$ implies $(\Delta^{s})$, i.e. any statistically $p$-ward continuous function is statistically ward continuous.
\begin{Thm} \label{TheopquasiCauchyimpliesquasiCauchy}   If $f$ is statistically $p$-ward continuous on a subset $A$ of $\mathbb{R}$ for some $p\in{\mathbb{N}}$, then it is statistically ward continuous on $A$.
\end{Thm}
\begin{pf}  If $p=1$, then it is obvious. So we would suppose that $p>1$. Take any statistically $p$-ward continuous function $f$ on $A$. Let $(\alpha_{k})$ be any statistical quasi-Cauchy sequence of points in $A$. Write $$(\xi_{i})=(\alpha_{1}, \alpha_{1},..., \alpha_{1}, \alpha_{2}, \alpha_{2},..., \alpha_{2},..., \alpha_{n}, \alpha_{n},..., \alpha_{n}, ...),$$ where the same term repeats $p$ times.
The sequence $$(\alpha_{1}, \alpha_{1},..., \alpha_{1}, \alpha_{2}, \alpha_{2},..., \alpha_{2},..., \alpha_{n}, \alpha_{n},..., \alpha_{n}, ...)$$ is also statistically quasi-Cauchy so it is statistically $p$-quasi-Cauchy. By the statistically $p$-ward continuity of $f$, the sequence
$$(f(\alpha_{1}), f(\alpha_{1}), ... ,f(\alpha_{1}),f(\alpha_{2}), f(\alpha_{2}),..., f(\alpha_{2}),..., f(\alpha_{n}), f(\alpha_{n}), ..., f(\alpha_{n}), ...)$$ is statistically $p$-quasi-Cauchy, where the same term repeats $p$-times.
Thus the sequence $$(f(\alpha_{1}), f(\alpha_{1}),..., f(\alpha_{1}), f(\alpha_{2}), f(\alpha_{2}),..., f(\alpha_{2}),..., f(\alpha_{n}), f(\alpha_{n}),..., f(\alpha_{n}), ...)$$ is also statistically $p$ quasi-Cauchy. It is easy to see that
$s_{t}-lim (f(\alpha_{n+p})-f(\alpha_{n}))=0$, which completes the proof of the theorem.

\end{pf}
\begin{Cor}   If $f$ is statistically $p$-ward continuous on a subset $A$ of $\mathbb{R}$, then it is continuous on $A$ in the ordinary case.
\end{Cor}
\begin{pf} The proof follows immediately from \cite[Theorem 3]{CakalliStatisticalwardcontinuity} so is omitted.
\end{pf}
\begin{Thm} Statistical $p$-ward continuous image of any statistically $p$-ward compact subset of $\mathbb{R}$ is statistically $p$-ward compact.
\end{Thm}
\begin{pf}
Let $f$ be a statistically $p$-ward continuous function, and $A$ be a statistically $p$-ward compact subset of $\mathbb{R}$. Take any sequence $\boldsymbol \beta=(\beta_{n})$ of terms in $f(E)$. Write $\beta_{n}=f(\alpha_{n})$ where $\alpha_{n}\in {E}$ for each $n \in{\mathbb{N}}$, $\boldsymbol\alpha=(\alpha_{n})$. statistically $p$-ward compactness of $A$ implies that there is a statistically $p$-quasi-Cauchy subsequence
$\boldsymbol \xi=(\xi_{k})=(\alpha_{n_{k}})$ of $\boldsymbol\alpha$. Since $f$ is statistically $p$-ward continuous, $(t_{k})=f(\boldsymbol \xi)=(f(\xi_{k}))$ is statistically $p$-quasi-Cauchy. Thus $(t_{k})$ is a statistically $p$-quasi-Cauchy subsequence of the sequence $f(\boldsymbol\alpha)$. This completes the proof of the theorem.
\end{pf}
\begin{Cor} Statistical $p$-ward continuous image of any $G$-sequentially connected subset of $\mathbb{R}$ is $G$-sequentially connected for a regular subsequential method $G$ (see \cite{CakalliOnGcontinuity}, \cite{MucukandCakalliGsequentiallyconnectednessfortopologicalgroupswithoperations}, and \cite{CakalliandMucukconnectednessviaasequentialmethod}).
\end{Cor}
\begin{pf} The proof follows from the preceding theorem, so is omitted (see \cite{CakalliSequentialdefinitionsofconnectedness} and \cite{MucukandSahanOnGsequentialcontinuity} for the definition of $G$-sequential connectedness and related concepts).
\end{pf}
\begin{Thm}  \label{TheouniformcontinuityimpliesthatfalpahansistatisticallypquasiCauchywheneveralphanispquasiCauchy} If $f$ is uniformly continuous on a subset $A$ of $\mathbb{R}$, then $(f(\alpha_{n}))$ is statistically $p$-quasi-Cauchy whenever $(\alpha_{n})$ is a $p$-quasi-Cauchy sequence of points in $A$.
\end{Thm}
\begin{pf}
Let $(\alpha_{n})$ be any $p$-quasi-Cauchy sequence of points in $A$. Take any $\varepsilon > 0$. Uniform continuity of $f$ on $A$ implies that there exists a $\delta >0$, depending on $\varepsilon$, such that $|f(x)-f(y)|< \varepsilon$ whenever  $|x-y|< \delta$ and $x, y\in{A}$. For this $\delta >0$, there exists an $N=N(\delta)$ such that $|\Delta_{p} \alpha_{n}|<\delta$ whenever $n>N$. Hence $|\Delta_{p} f(\alpha_{n})|<\varepsilon$ if $n>N$. Thus $\{k\leq n: |\Delta_{p} f(\alpha_{k})|\geq \varepsilon\}\subseteq \{1, 2, ..., N\}$. Therefore
$\lim_{n\rightarrow\infty}\frac{1}{n}|\{k\leq n:|\Delta_{p} f(\alpha_{k})|\geq \varepsilon\}|\leq \lim_{n\rightarrow\infty}\frac{1}{n}|\{k\leq N:k\in{\mathbb{N}}\}|=0$.
It follows from this that $(f(\alpha_{n}))$ is a statistically $p$-quasi-Cauchy sequence. This completes the proof of the theorem.
\end{pf}
\begin{Cor} If $f$ is slowly oscillating continuous on a bounded subset $A$ of $\mathbb{R}$, then $(f(\alpha_{n}))$ is statistically $p$-quasi-Cauchy whenever $(\alpha_{n})$ is a $p$ quasi-Cauchy sequence of points in $A$.
\end{Cor}
\begin{pf}
If $f$ is a slowly oscillating continuous function on a bounded subset $A$ of $\mathbb{R}$, then it is uniformly continuous on $A$ by \cite[Theorem 2.3]{CanakandDikNewtypesofcontinuities}. Hence the proof follows from Theorem \ref{TheouniformcontinuityimpliesthatfalpahansistatisticallypquasiCauchywheneveralphanispquasiCauchy}.
\end{pf}
It is well-known that any continuous function on a compact subset $A$ of $\mathbb{R}$ is uniformly continuous on $A$. We have an analogous theorem for a statistically $p$-ward continuous function defined on a statistically $p$-ward compact subset of $\mathbb{R}$.
\begin{Thm}
If a function is statistically $p$-ward continuous on a statistically $p$-ward compact subset of $\mathbb{R}$, then it is uniformly continuous on $A$.
\end{Thm}
\begin{pf}
Suppose that $f$ is not uniformly continuous on $A$ so that there exist an $\epsilon _{0}>0$ and sequences $(\alpha_{n})$ and $(\beta_{n})$ of points in $A$ such that $|\alpha_{n}-\beta_{n}|<1/n$ and $|f(\alpha_{n})-f(\beta_{n})|\geq \epsilon _{0}$ for all $n \in \mathbb{N}$. Since $A$ is statistically $p$-ward compact, there is a subsequence $(\alpha_{n_{k}})$ of $(\alpha_{n})$ that is statistically $p$-quasi-Cauchy. On the other hand, there is a subsequence $(\beta_{n_{k_{j}}})$ of $(\beta_{n_{k}})$ that is statistically $p$-quasi-Cauchy as well. It is clear that the corresponding sequence $(a_{n_{k_{j}}})$ is also statistically $p$-quasi-Cauchy, since \\
$\{j\leq n: |\alpha_{n_{k_{j+p}}}-\alpha_{n_{k_{j}}}|\geq{\varepsilon}\}\subseteq \{j\leq n: |\alpha_{n_{k_{j+p}}}-\beta_{n_{k_{j+p}}}|\geq{\frac{\varepsilon}{3}}\} \cup \{j\leq n: |\beta_{n_{k_{j+p}}}-\beta_{n_{k_{j}}}|\geq{\frac{\varepsilon}{3}}\}\cup \{j\leq n: |\beta_{n_{k_{j}}}-\alpha_{n_{k_{j}}}|\geq{\frac{\varepsilon}{3}}\}$\\
for every $n\in {\mathbb{N}}$, and for every $\varepsilon>0$.  Now the sequence \\$$(\alpha_{n_{k_{1}}},\alpha _{n_{k_{1}}},..., \alpha_{n_{k_{1}}}, \beta_{n_{k_{1}}}, \beta_{n_{k_{1}}},... ,\beta_{n_{k_{1}}},  ..., \alpha_{n_{k_{j}}}, \alpha_{n_{k_{j}}},...,\alpha_{n_{k_{j}}}, \beta_{n_{k_{j}}}, \beta_{n_{k_{j}}}, ..., \beta_{n_{k_{j}}},...)$$ is statistically $p$-quasi-Cauchy while the sequence
$$(f(\alpha_{n_{k_{1}}}),f(\alpha _{n_{k_{1}}}),...,f(\alpha_{n_{k_{1}}}),f(\beta_{n_{k_{1}}}),f(\beta_{n_{k_{1}}}),... ,f(\beta_{n_{k_{1}}}),  ...,f(\alpha_{n_{k_{j}}}),f(\alpha_{n_{k_{j}}}),...,f(\alpha_{n_{k_{j}}}), f(\beta_{n_{k_{j}}}),f(\beta_{n_{k_{j}}}),...,f(\beta_{n_{k_{j}}}),...)$$
is not statistically $p$-quasi-Cauchy where same term repeats $p$-times. Hence this establishes a contradiction, so completes the proof of the theorem.
\end{pf}
\begin{Cor} If a function defined on a bounded subset of $\mathbb{R}$ is statistically $p$-ward continuous, then it is uniformly continuous.
\end{Cor}
We note that when the domain of a function is restricted to a bounded subset of $\mathbb{R}$, statistically $p$-ward continuity implies not only ward continuity, but also slowly oscillating continuity.
\section{Conclusion}
In this paper, we introduce statistically $p$-quasi Cauchy sequences, and investigate conditions for a statistically $p$ ward continuous real function to be uniformly continuous, and prove some other results related to these kinds of continuities and some other kinds of continuities. It turns out that statistically $p$-ward continuity implies uniform continuity on a bounded subset of $\mathbb{R}$.  The results in this paper not only involves the related results in \cite{CakalliStatisticalwardcontinuity} as a special case for $p=1$, but also some interesting results which are also new for the special case $p=1$. The statistically $p$-quasi Cauchy concept for $p>1$ might find more interesting applications than statistical quasi Cauchy sequences to the cases when statistically quasi Cauchy does not apply. For a further study, we suggest to investigate statistically $p$-quasi-Cauchy sequences of soft points and statistically $p$-quasi-Cauchy sequences of fuzzy points. However due to the change in settings, the definitions and methods of proofs will not always be analogous to those of the present work (for example see \cite{ArasandSonmezandCakalliAnapproachtosoftfunctions}, \cite{CakalliandPratulFuzzycompactnessviasummability}, \cite{ErdemandArasandCakalliandSonmezSoftmatricesonsoftmultisetsinanoptimaldecisionprocess}, and  \cite{KocinacSelectionpropertiesinfuzzymetricspaces}). We also suggest to investigate statistically $p$-quasi-Cauchy double sequences of points in $\mathbb{R}$ (see \cite{PattersonandSavasAsymptoticequivalenceofdoublesequences}, \cite{{PattersonandCakalliQuasiCauchydoublesequences}}, \cite{DjurcicandKocinacandZizovicDoublesequencesandselections}, and \cite{CakalliandPattersonFunctionspreservingslowlyoscillatingdoublesequences} for the related definitions in the double case). For another further study, we suggest to investigate statistically $p$-quasi-Cauchy sequences in abstract metric spaces (see \cite{CakalliandSonmezandGencOnanequivalenceoftopologicalvectorspacevaluedconemetricspacesandmetricspaces}, \cite{PalandSavasandCakalliIconvergenceonconemetricspaces}, \cite{CakalliandSonmezSlowlyoscillatingcontinuityinabstractmetricspaces}, \cite{SonmezOnparacompactnessinconemetricspaces}, and  \cite{YayingandHazarikandCakalliNewresultsinquasiconemetricspaces}).

\end{document}